\documentstyle[11pt,leqno,rotate]{article}
\input amssym.def
\title{Upper bounds for the first eigenvalue of the Dirac operator on surfaces. \footnote {Supported by the SFB 288 of the DFG.}}
\author{Ilka Agricola and Thomas Friedrich \\
}
\date{}

\newcommand{\D}{\displaystyle}
\newcommand{\upsp}{\phantom{l}}
\newcommand{\downsp}{\phantom{q}}

\begin{document}

\begin{center}
\maketitle
{\small \it Humboldt-Universit\"at zu Berlin, Institut f\"ur Reine Mathematik,\\
Ziegelstra\ss e 13a, D-10099 Berlin, Germany\\
e-mail: agricola@mathematik.hu-berlin.de, friedric@mathematik.hu-berlin.de}
\end{center}

\mbox{} \hrulefill \mbox{}\\

{\bf Abstract}\\
In this paper we will prove new extrinsic upper bounds for the eigenvalues of the Dirac operator on  an isometrically immersed surface $M^2 \hookrightarrow {\Bbb R}^3$ as well as intrinsic bounds for 2-dimensional compact manifolds of genus zero and genus one. Moreover, we compare the different estimates of the eigenvalue of the Dirac operator for special families of metrics.\\

{\it Subj. Class.:} Differential geometry.

\mbox{} \hrulefill \mbox{}\\
\section{Introduction}

The Dirac operator $D$ acting on spinor fields defined over a 2-dimensional, compact, oriented Riemannian manifold $(M^2,g)$ with a fixed spin structure has a non-trivial kernel in general. Therefore, lower bounds for the eigenvalues of $D$ are not known in case the genus of $M^2$ is positive. The genus zero case is an exceptional one: using the uniformization theorem for simply-connected Riemann surfaces, we conclude that any metric $g$ on $S^2$ is conformally equivalent to the standard metric $g_o$ of $S^2$. Since the dimension of the space of all harmonic spinors depends on the conformal structure only,  it turns out that,  for any metric $g$ on $S^2$, there are no harmonic spinors. This observation yields a lower bound for the first eigenvalue $\lambda_1^2$ of $D^2$ proved by J. Lott (1986) and Chr. B\"ar (1992): the inequality

\[ \frac{4 \pi}{\mbox{vol} \,  (S^2, g)} \le \lambda_1^2 \]

holds for any Riemannian metric on $S^2$ (see [1], [11]).\\

On the other hand, several upper bounds for $\lambda_1^2$ depending on different geometric data are known. Intrinsic upper bounds involving the injectivity radius and the Gaussian curvature have been obtained by H. Baum (see [5]) and Chr. B\"ar (see [2]). In case the Riemannian surface $(M^2,g)$ is isometrically immersed into the 3-dimensional Euclidean space ${\Bbb R}^3$, one has extrinsic upper bounds depending on the ${\cal C}^0$-norm of the principal curvatures $\kappa_1, \kappa_2$ of the surface (see [5]). Denote by $H= (\kappa_1 + \kappa_2)/{2}$ the mean curvature. Then the following estimate for $ \lambda_1^2$ depending on the $L^2$-norm of the mean curvature $H$ is well-known (see [6], [3]):

\[ \lambda^2_1 \le \frac{\D \int_{M^2 \downsp}  H^2 dM^2}{\mbox{vol} \, (M^2,g)} . \]

In the present paper we will prove stronger extrinsic upper bounds for $\lambda_1^2$ in case of an isometrically immersed surface  $M^2 \hookrightarrow {\Bbb R}^3$ of arbritrary genus as well as an intrinsic upper bound for genus zero and genus one. Moreover, we will compare the different estimates of the eigenvalue of the Dirac operator for special families of metrics. \\

The extrinsic upper bound in case of a surface isometrically immersed into ${\Bbb R}^3$ depends on two smooth functions $f: M^2 \to {\Bbb R}$ and $G: {\Bbb R} \to {\Bbb R}$. \\

{\bf Theorem 1:} {\it The first eigenvalue $\lambda_1^2$ of the square of the Dirac operator on a surface $M^2 \hookrightarrow {\Bbb R}^3$ is bounded by}

\[ \lambda_1^2 \le \frac{\D \int_{M^2 \downsp} H^2 (f^2 + G^2 (f)) dM^2 + \int_{M^2} | \mbox{grad} \,  f|^2 (1+ [G' (f)]^2 )dM^2}{ \D \int_{M^2}^{\upsp}  ( f^2 + G^2 (f) ) dM^2} ,  \]

{\it where $f: M^2 \to {\Bbb R}$, $G: {\Bbb R} \to {\Bbb R}$ are smooth functions and $G'$ denotes the derivative of $G$.}\\

Suppose now that $(M^2,g)$ is a two-dimensional Riemannian manifold diffeomorphic to $S^2$. Denote by $g_o$ the standard metric of $S^2$. Then there exists a uniformization map, \mbox{i.e.,  a} conformal diffeomorphism $\Phi : S^2 \to M^2$. Let us introduce the function $h_{\Phi} : S^2 \to {\Bbb R}$ by the formula

\[ \Phi^* (g) = h_{\Phi}^4 g_o . \]

\newcommand{\dd}{\delta^{\mbox{\tiny Dir}}_c}

The set ${\cal U} (S^2, M^2)$ of all uniformization maps preserving the orientation can be parametrised by the elements  of the connected component of the group of all conformal diffeomorphisms of $S^2$, i.e.,  ${\cal U} (S^2, M^2) \approx SL (2, {\Bbb C})$. We introduce a new invariant $\dd (M^2,g)$ defined in a similar way as the conformal volume of a Riemann surface (see [10]):

\[ \dd (M^2,g) = \inf \left\{ \int\limits_{S^2} \frac{| \mbox{grad} (h_{\Phi})|^2}{h_{\Phi}^2} d S^2 : \Phi \in {\cal U} (S^2, M^2)  \right\} . \]

The vector field $\mbox{grad} (h_{\Phi})$ is the gradient of the function $h_{\Phi} :S^2 \to {\Bbb R}$ with respect to the standard metric of $S^2$.\\

{\bf Theorem 2:} {\it Let $(M^2,g)$ be a two-dimensional Riemannian manifold diffeomorphic to the sphere $S^2$. Then}


\[ 0 \le \, \, \lambda_1^2 - \frac{4 \pi}{\mbox{vol} \, (M^2,g)} \, \,  \le  \, \, \frac{\dd (M^2,g)}{\mbox{vol} \, (M^2,g)} \]

{\it holds.}\\

The same method applies to Riemannian metrics on the two-dimensional torus $T^2$. The spin structures of $T^2$ are described by pairs $(\varepsilon_1, \varepsilon_2)$ of numbers $\varepsilon_i =0,1$, the trivial spin structure corresponding to the pair $(\varepsilon_1 , \varepsilon_2) = (0,0)$. Let $\Gamma$ be a lattice in ${\Bbb R}^2$ with basis $v_1, v_2$ and denote by $v_1^*, v_2^*$ the dual basis of the dual lattice $\Gamma^*$. We will compare the flat metric $g_o$ on the torus $T^2 = {\Bbb R}^2/ \Gamma$ with a conformally equivalent metric $g = h^4 g_o$.\\

{\bf Theorem 3:} {\it Let $(M^2,g)$ be a two-dimensional Riemannian manifold conformally equivalent to the flat torus $T^2$ and equipped with the trivial spin structure. Then the Dirac operator on $(M^2,g)$ has a two-dimensional kernel. Moreover, the first positive eigenvalue $\lambda_1^2 (g)$ of $D^2$ on $(M^2,g)$ is bounded by }

\[ \lambda_1^2 (g) \le \frac{\D \int_{T^2 \downsp}  \left\{ \lambda_1^2 (g_o) + \frac{4}{h^2} | \mbox{grad} \, (h)|^2 \right\} \frac{1}{h^6} dT^2}{\D \int_{T^2}^{\upsp} \frac{1}{h^2} dT^2} . \]

{\bf Theorem 4:} {\it Let $(M^2,g)$ be a two-dimensional Riemannian manifold conformally equivalent to the flat torus $T^2$. In case the spin structure $(\varepsilon_1 , \varepsilon_2) \not= (0,0)$ is non-trivial,  the Dirac operator has a trivial kernel and $\lambda_1^2 (D)$ is bounded by

\[ \lambda_1^2 (D) \le \pi^2 \, | \varepsilon_1 v^*_1 + \varepsilon_2 v_2^* |^2 \, \frac{\D \int_{T^2 {\downsp}}  \frac{1}{h^2} dT^2}{\D \int_{T^2}^{\upsp}  h^2 dT^2} \, \, . \]

Moreover, the inequality}

\[ \lambda_1^2 (D) \mbox{vol} \, (M^2,g) \le \lambda_1^2 (g_o) \mbox{vol} \, (T^2,g_o) + \int_{T^2} \frac{|\mbox{grad} \, (h)|^2}{h^2} dT^2 \]

{\it with} 

\[ \lambda_1^2 (g_o) \mbox{vol} \, (T^2, g_o) = \pi^2 \frac{| \varepsilon_1 v_1^* + \varepsilon_2 v_2^* |^2}{\sqrt{|v_1^*|^2 |v_2^*|^2 - \langle v_1^*, v_2^* \rangle}} \]

{\it holds.}\\

We shall apply the previous results to two families of surfaces of special interest. Let us first consider  the ellipsoid

\[ E(a) = \left\{ (x,y,z) \in {\Bbb R}^3: x^2 +y^2 + \frac{z^2}{a^2} =1 \right\}  \, . \]

A calculation of the volume yields that the lower bound ${4 \pi}/{\mbox{vol} \, (E(a))}$ for $\lambda_1^2 (a)$ is a monotone decreasing function of the parameter $a$:

\[ \lim\limits_{a \to 0} \frac{4 \pi}{\mbox{vol} \, (E(a))} = 2 \quad , \quad \lim\limits_{a \to \infty} \frac{4 \pi}{\mbox{vol} \, (E(a))} = 0 . \]

Using the upper bounds for $\lambda_1^2 (a)$ already known, we cannot control the behaviour of $\lambda_1^2 (a)$ for small or large values of the parameter $a$. For example, the $L^2$-bound given by the mean curvature $H$ has the following limits: 

\[  \lim_{a \to 0} \frac{\D \int_{E(a) \downsp} H^2 dE(a)}{\mbox{vol} \, (E(a))} = \infty \quad , \quad   \lim_{a \to 0} \frac{\D \int_{E(a) } H^2 dE(a)}{\mbox{vol} \, (E(a))} = \frac{1}{2} . \]

Now, a combination of our stronger extrinsic and intrinsic upper bounds for the first eigenvalue of the Dirac operator yields the following  improvement for the ellipsoid:\\

{\bf Theorem 5:} {\it The first eigenvalue $\lambda_1^2$ of $D^2$ on the ellipsoid $E(a)$ satisfies 
\begin{itemize}
\item[1.)]  $2 \le \overline{\lim\limits_{a \to 0}} \, \lambda_1^2 (a) \le \frac{3}{2} + \ln 2 \approx 2,2$;
\item[2.)]  $\lim\limits_{a \to \infty} \lambda_1^2 (a) =0$;
\item[3.)]  $\lambda_1^2 (a) \stackrel{<}{\sim} \frac{2 \ln (2) + 3}{\pi} \frac{1}{a} $ \, \, for $a \to \infty.$
\end{itemize}}

In the last part of this paper we apply our estimates to a tube of radius $r$ around a circle of curvature $\kappa$, i.e., a ''round'' torus. Parametrizing the spin structure as before, the inequalities for $\lambda_1^2 (\kappa, r)$ allow us to prove, in particular,

\[ \lim\limits_{r \to 0} \lambda_1^2 (\kappa,r) \mbox{vol} \, (\kappa,r) = \lim\limits_{\kappa \to 0} \lambda_1^2 ( \kappa,r) \mbox{vol} \, (\kappa,r)=0 \]

for the spin structure $(\varepsilon_1, \varepsilon_2)=(1,0)$ and 

\[ \overline{\lim\limits_{r \kappa \to 1}} \lambda_1^2 (\kappa,r) \mbox{vol} \, (\kappa, r) \le \pi^2 \]

for the spin structure $(\varepsilon_1, \varepsilon_2)=(0,1)$ (for these two spin structures, no upper bounds were available before). However, they turn out to yield no improvement for the induced spin structure $(\varepsilon_1, \varepsilon_2)=(1,1)$; thus, in this case, the classical bound involving the integral over $H^2$ divided by the volume is still the best one available. \\

\section{Extrinsic upper bounds}

Let $M^2$ be a compact, oriented surface isometrically immersed into the Euclidean space ${\Bbb R}^3$ and denote by $\vec{N} (m)$ the unit normal vector of $M^2$ at the point $m \in M^2$. The restriction $\Phi_{|M^2}$ of a spinor field $\Phi$ defined on ${\Bbb R}^3$ is a spinor field on the surface $M^2$. Let $\Phi$ be a parallel spinor on ${\Bbb R}^3$. Then the spinor field

\[ \varphi^* = \frac{1}{2} (1-i) \Phi_{|M^2} + \frac{1}{2} (-1+i) \vec{N} \cdot \Phi_{|M^2} \]

is of constant length on $M^2$ and satisfies  the two-dimensional Dirac equation

\[ D( \varphi^*)= H \varphi^*  , \]

where $H$ denotes the mean curvature of the surface (see [9]). Thus, starting with two parallel spinors $\Phi_1, \Phi_2$ with 

\[ | \Phi_1 |= | \Phi_2 |=1 \quad \mbox{and} \quad \langle \Phi_1 , \Phi_2 \rangle =0  \, , \]

we obtain two solutions $\varphi_1^*, \varphi^*_2$ of the Dirac equation

\[ D(\varphi_{\alpha}^*)= H \varphi_{\alpha}^* \quad , \quad \alpha =1,2 \]

such that $| \varphi_1^* (m)| = | \varphi_2^* (m)| =1$ and $\langle \varphi_1^* (m) , \varphi_2^* (m) \rangle =0$ holds at any point \linebreak $m \in M^2$. Given two real-valued functions $f,g : M^2 \to {\Bbb R}$ we consider the spinor field

\[ \psi = f \varphi_1^* + g \varphi_2^*  . \]

After applying  the Dirac operator to $\psi$

\[ D( \psi)= H \psi + \mbox{grad} \, (f) \cdot \varphi_1^* + \mbox{grad} \, (g) \cdot \varphi_2^*  , \]

a direct calculation yields the formula

\[ |D(\psi)|^2 = H^2 (f^2 + g^2 ) + | \mbox{grad} \, (f)|^2 + |\mbox{grad} \, (g)|^2 - 2 \, \mbox{Re} \,  (\mbox{grad} \, (f) \cdot \mbox{grad} \, (g) \cdot \varphi_2^* , \varphi_1^*) . \]

In case the vector fields $\mbox{grad} \, (g)$ and $\mbox{grad} \, (f)$ are parallel,  the last term in this formula vanishes since $\varphi_1^*$ and $\varphi_2^*$ are orthogonal. In this case the Rayleigh quotient coincides with

\[ \frac{\D \int_{M^2 \downsp} |D(\psi)|^2}{\D \int_{M^2}^{\upsp} |\psi|^2} =  \frac{ \D \int_{M^2 \downsp}  H^2 (f^2 + g^2) dM^2 + {\D \int_{M^2}} \left( |\mbox{grad} \, (f)|^2 + |\mbox{grad} \, (g)|^2 \right) dM^2}{\D \int_{M^2}^{\upsp} (f^2 + g^2) dM^2} . \]

The condition for the gradients of the functions $f$ and $g$ is satisfied for example if  $g$ is a function depending on $f$, i.e.,  $g=G(f)$. Finally, we have  proved Theorem 1.\\

\section{Intrinsic upper bounds for a surface diffeomorphic to \boldmath $ S^2  \, \mbox{or} \, T^2$ \unboldmath}

Let $(M^2, g_o)$ be a compact, oriented 2-dimensional Riemannian spin manifold and denote by $D_o$ its Dirac operator. Moreover,  consider a conformally equivalent metric

\[ {g} = h^4 g_o . \]

The corresponding Dirac operator ${D}$ is related with $D_o$ by the formula (see [4]) 

\[ {D} = \frac{1}{h^2} D_o + \frac{\mbox{grad} \, (h)}{h^3}  . \]

Consequently, the equation ${D} (\psi) = {\lambda} \psi$ is equivalent to

\[ D_o(\psi)= {\lambda} h^2 \psi - \frac{1}{h} \mbox{grad} (h) \cdot \psi . \]

For any spinor field $\psi$ we compute the $L^2$-norm of ${D} (\psi)$: 

\[ \int\limits_{M^2} | {D} (\psi) |^2 d {M}^2 = \int\limits_{M^2} \left\{ |D_o(\psi) |^2 + \frac{| \mbox{grad} \, (h)|^2}{h^2} |\psi|^2 + \frac{2}{h} \mbox{Re} \,  (\mbox{grad} \, (h) \cdot \psi, D_o( \psi)) \right\} dM_o^2 . \]

Suppose now that $\psi$ is an eigenspinor of the Dirac operator $D_o$ with eigenvalue $\lambda_i$. Then $\mbox{Re} \,  (\mbox{grad} \, (h) \cdot \psi, D_o(\psi))=0$ and we obtain the formula

\[ \int\limits_{M^2} | {D} (\psi) |^2 d {M}^2 = \int\limits_{M^2}  \left\{ \lambda_i^2 + \frac{|\mbox{grad} \, (h)|^2}{h^2} \right\} |\psi|^2 dM_o^2 . \]

Hence, the first eigenvalue $\lambda_1^2 ({D})$ of the Dirac operator is bounded by

\[ \lambda_1^2 ({D}) \le \inf\limits_{\lambda_i} \, \, \, \inf\limits_{D_o (\psi)= \lambda_i \psi} \frac{\D \int_{M^2} \downsp \left\{ \lambda_i^2 + \frac{|\mbox{grad} \, (h)|^2}{h^2} \right\} |\psi|^2 dM^2_o}{\D \int_{M^2}^{\upsp} |\psi|^2 h^4 dM^2_o} . \]

Let us now discuss the special case that $(M^2,g_o)$ is the two-dimensional sphere with its standard metric and ${g}$  a conformally equivalent metric. The first eigenvalue of the Dirac operator on $S^2$ is $\lambda_1 =1$. Moreover, the corresponding eigenspinor $\psi$ is a real Killing spinor satisfying the differential equation

\[ \nabla_X (\psi) = - \frac{1}{2} X \cdot \psi \quad , \quad X \in T(S^2) . \]

In particular, the length of $\psi$ is constant and we obtain the inequality

\[ \lambda_1^2 ({D}) \le \frac{4 \pi}{\mbox{vol} \, (S^2,g)} + \frac{\D \int_{S^2 \downsp} \frac{| \mbox{grad} \, (h)|^2}{h^2} dS^2}{\mbox{vol} \, (S^2,g)} . \]

Starting with a surface $(M^2,g)$ diffeomorphic to $S^2$, the latter inequality holds for any uniformization, i.e.,  for any conformal diffeomorphism $\Phi : S^2 \to M^2$ such that $\Phi^* (g) = h_{\Phi}^4 g_o$. In particular, we have  proved Theorem 2.\\

{\bf Remark:} For any conformal diffeomorphism $\psi \in SL (2, {\Bbb C})$ of the two-dimensional sphere $S^2$ we denote by $h_{\psi} : S^2 \to {\Bbb R}$ the function defined by the equation

\[ \psi^* (g_o)= h_{\psi}^4 g_o . \]

Let $f:S^2 \to {\Bbb R}$ be a smooth function. Then we define the number

\[ \delta_c^{\tiny \mbox{Dir}} (f) = \inf \left\{  \int_{S^2} | \mbox{grad} (f \circ \psi) + \mbox{grad} (\log (h_{\psi}))|^2 dS^2: \psi \in SL (2, {\Bbb C}) \right\} \, \, . \]

In case of a uniformization $\Phi: S^2 \to M^2$ such that $\Phi^* (g) = h^4_{\Phi} g_o$, we have

\[ \int_{S^2} \frac{| \mbox{grad} (h_{\Phi})|^2}{h_{\Phi}^2} dS^2 = \int_{S^2} |€\mbox{grad} (\log (h_{\Phi}))|^2 dS^2 \]

and, consequently, for the quantity $\delta^{\tiny \mbox{Dir}}_c (M^2,g)$ defined in the introduction,  the relation 

\[ \delta^{\tiny \mbox{Dir}}_c (M^2,g) = \delta^{\tiny \mbox{Dir}}_c (\log (h_{\Phi})) . \]

We consider  the case that $(M^2,g)$ is the flat  torus $T^2 =({\Bbb R}^2 /\Gamma, g_o)$ given by a lattice $\Gamma$ in ${\Bbb R}^2$ with trivial spin structure. In this case  there are two parallel spinor fields $\varphi^+$ and $\varphi^-$ of constant length and the first non-trivial eigenvalue $\lambda_1^2 (g_o)$ of the square of the Dirac $D_o$ operator on $T^2$ is

\[ \lambda_1^2 (g_o) = 4 \pi^2 \min \left\{ |v^*|^2: \, \, \, 0 \not= v^* \in \Gamma \right\}  , \]

where $\Gamma^*$ denotes the dual lattice (see [7]). Suppose now that $g$ is a metric on $M^2$ conformally equivalent to $g_o$, $g=h^4 g_o$. Then the kernel of the corresponding Dirac operator is again two-dimensional and spanned by the spinor fields $\frac{1}{h} \varphi^+, \frac{1}{h} \varphi^-$. Fix a spinor field $\psi$ such  that $D_o (\psi)= \lambda_1 (g_o) \psi$.  Then the length of $\psi$ is constant, i.e., $|\psi| \equiv 1$. The spinor field $\psi^* = \psi / h^3$ is orthogonal to the kernel of the Dirac operator $D$ with respect to the $L^2$-norm of the metric $g$. Indeed, we have

\[ \int_{M^2} \left( \psi^*, \frac{1}{h} \varphi^{\pm} \right) dM^2 = \int_{T^2} \left( \frac{1}{h^3} \psi, \frac{1}{h} \varphi^{\pm} \right) h^4 dT^2 = \]

\[ = \frac{1}{\lambda_1 (g_o)} \int_{T^2} \left( D_o (\psi), \varphi^{\pm} \right) dT^2 = \frac{1}{\lambda_1 (g_o)} \int_{T^2} \left( \psi, D_o (\varphi^{\pm}) \right) dT^2 = 0 . \]

This observation yields the inequality

\[ \lambda_1^2 (g) \le \frac{\D \int_{M^2 \downsp} |D(\psi^*)|^2 dM^2}{\D \int_{M^2}^{\upsp} |\psi^*|^2 dM^2} \]

for the first non-trivial eigenvalue of $D^2$ on $(M^2,g)$. Moreover, we have

\[ \int_{M^2} |\psi^*|^2 dM^2 = \int_{T^2} \frac{1}{h^6} h^4 dT^2 = \int_{T^2} \frac{1}{h^2} dT^2 \]

and 

\[ \int_{M^2} |D(\psi^*)|^2 dM^2 = \int_{T^2} \left\{ |D_o (\psi^*)|^2 + \frac{\mbox{grad} \, (h)}{h^2} |\psi^*|^2+ \frac{2}{h} \mbox{Re} \, \left( \mbox{grad} \, (h) \psi^*, D_o(\psi^*) \right) \right\} dT^2 . \]

Since the equation

\[ D_o (h^3 \psi^*)= D_o (\psi) = \lambda_1 (g_o) \psi = \lambda_1 (g_o) h^3 \psi^* \]

can be rewritten in the form

\[ D_o (\psi^*)= \lambda_1 (g_o) \psi^* - \frac{3}{h} \, \, \mbox{grad} \, (h) \psi^* , \]

we obtain the formulas

\[ \frac{2}{h} \mbox{Re} \, (\mbox{grad} \, (h) \psi^*, D_o (\psi^*))= - \frac{6}{h^2} \, \, |\mbox{grad} \, (h) |^2 |\psi^*| \]

and

\[ |D_o (\psi^*)|^2 = \left\{ \lambda_1^2 (g_o) + \frac{9}{h^2} \, \, |\mbox{grad} \, |^2 \right\} |\psi^*|^2 . \]

Altogether, this implies

\[
\int_{M^2} |D(\psi^*)|^2 dM^2 =  \int_{T^2} \left\{ \lambda_1^2 (g_o) + \frac{4}{h^2} \, \, |\mbox{grad} (h)|^2 \right\} \frac{1}{h^6} dT^2
\]

and it proves Theorem 3, in particular.\\

Let us now consider the case that the spin structure on $(M^2,g) \approx T^2$ is non-trivial. Then the Dirac operator has no kernel and the eigenspinors of the Dirac operator $D_o$ on $T^2$ are again of constant length (see [7]). Then our method provides the inequality

\[ \lambda_1^2 (g) \le \frac{\lambda_1^2 (g_o) \mbox{vol} \, (T^2, g_o)}{\mbox{vol} \, (M^2,g)} + \frac{\D \int_{T^2 \downsp} \frac{|\mbox{grad} (h)|^2}{h^2} \, \, dT^2}{\mbox{vol} \, (M^2,g)} . \]

The Gaussian curvature $G$ of the metric $g$ is given by

\[ h^4 G = - 2 \Delta (\log (h))  , \]

where $\Delta$ denotes the Laplacian with respect to the flat metric. We integrate this latter equation:

\[ \int_{M^2} G \cdot \log (h) dM^2 = \int_{T^2} h^4 G \cdot \log (h) dT^2 = - 2 \int_{T^2} \frac{|\mbox{grad} (h)|^2}{h^2} \, \, dT^2  \]

thus obtaining 

\[ \lambda_1^2 (g) \le \frac{\lambda_1^2 (g_o) \mbox{vol} \, (T^2, g_o)}{\mbox{vol} \, (M^2,g)} - \frac{1}{2}  \, \, \frac{\D \int_{M^2 \downsp} G \cdot \log (h) dM^2}{\mbox{vol} \, (M^2,g)}  , \]

where $G$ denotes the Gaussian curvature of $(M^2,g)$.\\

However, we can use a more delicate comparison for the Dirac operator depending on the spin structure $(\varepsilon_1, \varepsilon_2)$. Consider the dual lattice $\Gamma^*$ with basis $v_1^*, v_2^*$ as well as the 1-form

\[ \omega = \pi i (dx, dy) \cdot (\varepsilon_1 v_1^* + \varepsilon_2 v_2^*) . \]

The Dirac operator $D^{(\varepsilon_1, \varepsilon_2)}$ corresponding to the spin structure $(\varepsilon_1, \varepsilon_2)$ on $(M^2,g)$ is related to the Dirac operator $D$ for the trivial spin structure by

\[ D^{(\varepsilon_1, \varepsilon_2)} = D+it , \]

where the vector field $t$ is dual with respect to the metric $g$ to the $1$-form $\omega$ (see [7]). Let $\varphi^+$ be the parallel spinor field with respect to the flat metric. Then $\psi = \frac{1}{h} \varphi^+$ is a harmonic spinor on $(M^2,g)$, i.e.,  $D(\psi) = 0.$ Therefore,  we obtain

\[ |D^{(\varepsilon_1, \varepsilon_2)} (\psi)|^2 = |t|^2_g |\psi|^2 = |\omega|^2_g |\psi|^2 . \]

In dimension $n=2$ the $L^2$-length of a  1-form depends only the conformal structure, i.e.,  if the metrics $g=h^4 g_o$ and $g_o$ are conformally equivalent, then for any 1-form $\omega$ the formula

\[ |\omega|^2_g dM^2_g =|\omega|^2_{g_o} dM^2_{g_o} \]

holds. Now we integrate:

\[
\int_{M^2} |D^{(\varepsilon_1, \varepsilon_2)} (\psi) |^2 dM^2 =  \int_{T^2} \frac{1}{h^2} |\omega|^2_{g_o} dT^2 = \pi^2 |\varepsilon_1 v_1^* + \varepsilon_2 v_2^* |^2 \, \, \int_{T^2} \frac{1}{h^2} \, d T^2 . 
\]

On the other hand, we have

\[ \int_{M^2} |\psi|^2 dM^2 = \int_{T^2} \frac{1}{h^2} h^4 dT^2 = \int_{T^2} h^2 dT^2  ; \]

finally, we obtain

\[ \frac{\D \int_{M^2 \downsp} |D^{(\varepsilon_1, \varepsilon_2)} (\psi) |^2 dM^2}{\D \int_{M^2}^{\upsp} |\psi|^2 dM^2} = \pi^2 |\varepsilon_1 v_1^* + \varepsilon_2 v_2^*|^2 \, \, \frac{\D \int_{T^2 \downsp} \frac{1}{h^2} dT^2}{\D \int_{T^2}^{\upsp} h^2 dT^2} . \]

This equality finishes the proof of Theorem 4.\\

\section{The first eigenvalue of the Dirac operator on the ellipsoid with \boldmath $S^1$\unboldmath-symmetry}

We now discuss the first eigenvalue of the Dirac operator on the ellipsoid $E(a) \subset {\Bbb R}^3$ with $S^1$-symmetry defined by the equation

\[ x^2 + y^2 + \frac{z^2}{a^2} = 1 . \]

For the calculations we will use the following convenient parametrization of $E(a)$:

\[ x = \sqrt{1- w^2} \cos \varphi \quad , \quad y = \sqrt{1-w^2} \sin \varphi \quad , \quad z= a \cdot w \quad ,  \]

where the parameters $(w, \varphi)$ are restricted to the intervals $-1 \le w \le1, \, \, 0 \le \varphi \le 2 \pi$. For brevity we introduce the function 

\[ \Delta_a (w) = (1-a^2) w^2 + a^2 . \]

Then the Riemannian metric $ds^2_a$, the Gaussian curvature $G$, the mean curvature $H$ and the volume form $dE(a)$ are given by the formulas:

\begin{itemize}
\item[1.)] \mbox{}  $ds_a^2 = \frac{\Delta_a (w)}{1-w^2} dw^2 + (1-w^2)d \varphi^2$; 
\item[2.)] \mbox{}  $H^2 = \frac{a^2}{4} \Delta^{-3}_a (w) \{ \Delta_a (w) +1 \}^2$;
\item[3.)] \mbox{}  $G= a^2 \Delta_a^{-2} (w)$;
\item[4.)] \mbox{}  $dE(a) = \Delta_a^{1/2} (w) dw \wedge d \varphi$.
\end{itemize}

\subsection{Evaluation of the extrinsic upper bounds}

We shall  use the extrinsic upper bound for the eigenvalue of the Dirac operator for the family of functions  $f_{\beta}$ defined by

\[ f_{\beta} = \Delta_a^{\beta} (w) \quad , \quad \beta > \frac{1}{2} . \]

Notice that $f_{\beta}$ is just the $\beta$-th power of (a multiple of) $1/\sqrt{G}$. The length of the gradient of the function $f_{\beta}$ on the ellipsoid is given by

\begin{itemize}
\item[5.)] \mbox{} $|€\mbox{grad} \, (f_{\beta})|^2 = 4 \beta^2 (1- a^2)^2 \Delta^{2 \beta - 3}_a (w) w^2 (1-w^2)$.
\end{itemize}

Let us first discuss the case that the parameter $a <1$ is small. Then $a^2 \le \Delta_a (w) \le 1$ holds and we can estimate the first  integral appearing in Theorem 1

\[ 0 \le \int_{E(a)} H^2 f_{\beta}^2 dE(a) \le 4 \pi a^2 \int\limits^1_0 \Delta_a^{2 \beta - 5/2} (w) dw . \]

The latter integral may be rewritten using the transformation $\sqrt{1-a^2} \, \, w= ax$,  thus yielding 

\[ 0 \le \int_{E(a)} H^2 f_{\beta}^2 dE(a) \le 4 \pi \frac{a^{4 \beta -2}}{\sqrt{1-a^2}} \int\limits^{\frac{1}{a}\sqrt{1-a^2}}_0 (1+x^2)^{2 \beta - 5/2}  dx . \]

We shall prove that for all $\beta > \frac{1}{2}$

\[ \lim\limits_{a \to 0} \int_{E(a)} H^2 f_{\beta}^2 dE(a)=0  . \]

Indeed, in case $\beta \ge \frac{5}{4}$,   we have $\Delta_a^{2 \beta - 5/2} (w) \le 1$ and the result follows immediately. If $\frac{3}{4} < \beta \le \frac{5}{4}$,  we use the inequality $a^2 \le \Delta_a (w)$, i.e.,  $\Delta^{2 \beta - 5/2}_a (w) \le a^{4 \beta - 5}$. Finally, consider the case that $\frac{1}{2} < \beta \le \frac{3}{4}$. Then one has $1 \le \frac{5}{2} - 2 \beta < \frac{3}{2}$ and, hence,  \linebreak $(1+x^2) \le (1+x^2)^{5/2 - 2 \beta}$, which implies

\[ \int\limits^{\frac{1}{a}\sqrt{1-a^2}}_0 (1+x^2)^{2 \beta - 5/2} dx \le \int\limits^{\infty}_0 \frac{dx}{1+x^2} < \infty  \]

and finishes the argument. In a similar way we show

\[ \lim\limits_{a \to 0} \int_{E(a)} f^2_{\beta} dE(a) = \lim\limits_{a \to 0} 4 \pi \int\limits^1_0 \Delta_a^{2 \beta + 1/2} (w) dw = 4 \pi \int\limits^1_0 w^{4 \beta +1} dw = \frac{2 \pi}{2 \beta +1} . \]

Finally,  we investigate the integrals

\[ \int_{E(a)} |\mbox{grad} (f_{\beta})|^2 dE(a)= 16 \pi (1-a^2)^2 \beta^2 \int\limits^1_0 \Delta_a^{2 \beta - 5/2} (w) w^2 (1-w^2) dw . \]

Using the Lebesgue theorem $(\beta > \frac{1}{2})$ we conclude

\[ \lim\limits_{a \to 0} \int_{E(a)} |\mbox{grad} (f_{\beta})|^2 dE(a) = 16 \pi \beta^2 \int\limits^1_0 w^{4 \beta - 3} (1-w^2) dw = 4 \pi \frac{\beta}{2 \beta - 1} . \]

Since the first eigenvalue $\lambda_1^2 (a)$ of the square of the Dirac operator on $E(a)$ is bounded by the expression 

\[ \lambda_1^2 (a) \le \frac{\D \int_{E(a)} \downsp H^2 f^2_{\beta} dE(a) + \int_{E(a)} |\mbox{grad} (f_{\beta})|^2 dE(a)}{\D \int_{E(a)}^{\upsp}  f^2_{\beta} \, dE(a)}  , \]

we obtain 

\[ \overline{\lim\limits_{a \to 0}} \lambda^2_1 (a) \le \frac{4 \pi \beta \cdot (2 \beta +1)}{(2 \beta - 1) 2 \pi} = 2 \frac{\beta (2 \beta +1)}{2 \beta - 1}  \]

in the limit $a \to 0$. The latter inequality holds for any $\beta > \frac{1}{2}$. For $\beta =1$ we obtain, for example,  the inequality

\[ \overline{\lim\limits_{a \to 0}} \lambda_1^2 (a) \le 6 \]

and the optimal parameter $\beta = \frac{1}{2} + \frac{1}{\sqrt{2}}$ yields the estimate

\[ \overline{\lim\limits_{a \to 0}} \lambda_1^2 (a) \le 3 + 2 \sqrt{2} \approx 5,8 . \]

Later, this result will be sharpened with the aid of the intrinsic bounds; however, we already get as a partial result that $\lambda^2_1$ remains bounded. \\
We now discuss the case of a large parameter $a$ \, $(a >1)$. It is convenient to  write $\Delta_a (w)$ in the form  $\Delta_a (w) =(a^2-1) \left[ \frac{a^2}{a^2-1} - w^2 \right]$. The  formulas 1.) - 5.)  used before imply

\[ \frac{\D \int_{E(a)} \downsp |\mbox{grad} (f_{\beta})|^2 dE(a)}{\D \int_{E(a)}^{\upsp}  f^2_{\beta} dE(a)} = 4 \beta^2 \frac{1}{a^2 - 1}  \cdot \frac{ \D \int\limits^1_0 \left[\frac{a^2}{a^2-1} - w^2 \right]^{2 \beta - 5/2} w^2 (1-w^2)dw}{\D \int\limits^1_0 \left[\frac{a^2}{a^2-1} - w^2 \right]^{2 \beta + 1/2} dw} . \]

We compute again its limit for $a \to \infty$:

\[ \lim\limits_{a \to \infty} \frac{\D \int_{E(a)} \downsp |\mbox{grad} (f_{\beta})|^2 dE(a)}{\D \int_{E(a)}^{\upsp}  f^2_{\beta} dE(a)} =0 . \]

Thus, the asymptotic behaviour is dominated by the second term of the estimate: 

\[ \lim\limits_{a \to \infty} \, \,  \frac{\D \int_{E(a)} \downsp H^2 f^2_{\beta} dE(a)}{\D \int_{E(a)}^{\upsp}  f^2_{\beta} dE(a)} = \frac{1}{4} \, \frac{\D \int\limits^1_0 [1-w^2]^{2 \beta - 1/2} dw}{\D \int\limits^1_0 [1-w^2]^{2 \beta + 1/2} dw} . \]

This yields  the inequality

\[ \overline{\lim\limits_{a \to \infty}} \, \, \lambda_1^2 (a) \le \frac{1}{4} \, \, \frac{\D \int\limits^1_0 [1-w^2]^{2 \beta - 1/2} dw}{\D \int\limits^1_0 [1-w^2]^{2 \beta + 1/2} dw} \]

for any $\beta > \frac{1}{2}$. The special value in case of the parameter $\beta =1$ can easily be calculated to be 

\[ \overline{\lim\limits_{a \to \infty}} \lambda_1^2 (a) \le \frac{3}{10} . \]

However, the inequality holds for any $\beta > \frac{1}{2}$; for  $\beta \to \infty$ we obtain the optimal result 

\[ \overline{\lim\limits_{a \to \infty}} \lambda^2_1 (a) \le \frac{1}{4} . \]

{\bf Remark:} Let us point out that, for $\beta=1$, the integral approximation of $\lambda_1^2 (a)$ is,  on both sides $a \to 0, \infty$,  not the best one among the extrinsic upper bounds considered, but we may come  very close to the optimal value using the family of functions $f_{\beta}$. The exact formula holding  for all parameters $0 <a< \infty$ is  in this case:

\[ \lambda_1^2 (a) \le \frac{\left( 2 + \frac{13}{8} a^2 + \frac{3}{16} a^4 \right) + \left( \frac{7}{2} a^2 - \frac{3}{2} a^4 - \frac{3}{16} a^6 \right) f(a)}{\left( \frac{1}{3} + \frac{5}{12} a^2 + \frac{5}{8} a^4 \right) - \frac{5}{8} a^6 f(a)} \]

where the function $f(a)$ is given by

\[ f(a) = \left\{ \begin{array}{ll} \frac{1}{\sqrt{1-a^2}} \ln \left( \frac{1 - \sqrt{1-a^2}}{a} \right) & a < 1\\ \mbox{}\\ - \frac{1}{\sqrt{a^2 - 1}} \arcsin \left( \frac{\sqrt{a^2-1}}{a} \right) & a > 1 \end{array} \right.  . \]

Figure  1 $(a \in [0,1[)$ and figure 2 $(a \in ]1, \infty [)$ give an overview of the different extrinsic bounds. The lower solid line is the only known lower bound proportional to the inverse of the volume due to Lott and B\"ar; the upper solid line is the well known upper bound involving  the integral over $H^2$ divided by the volume. The short dashed curve corresponds to $\beta = 1/2$ in our family of functions; as seen before, this is the maximal value for $\beta$ for which the curve does not remain bounded as $a \to 0$. Its limit for $a \to \infty$ is $1/3$. Finally, the long dashed curve is the upper bound for $\beta =1$ as discussed previously.\\

\input epsfig.sty

\vspace{-1cm}
\begin{center}
\[
\begin{rotate}[r]{\epsfig{figure=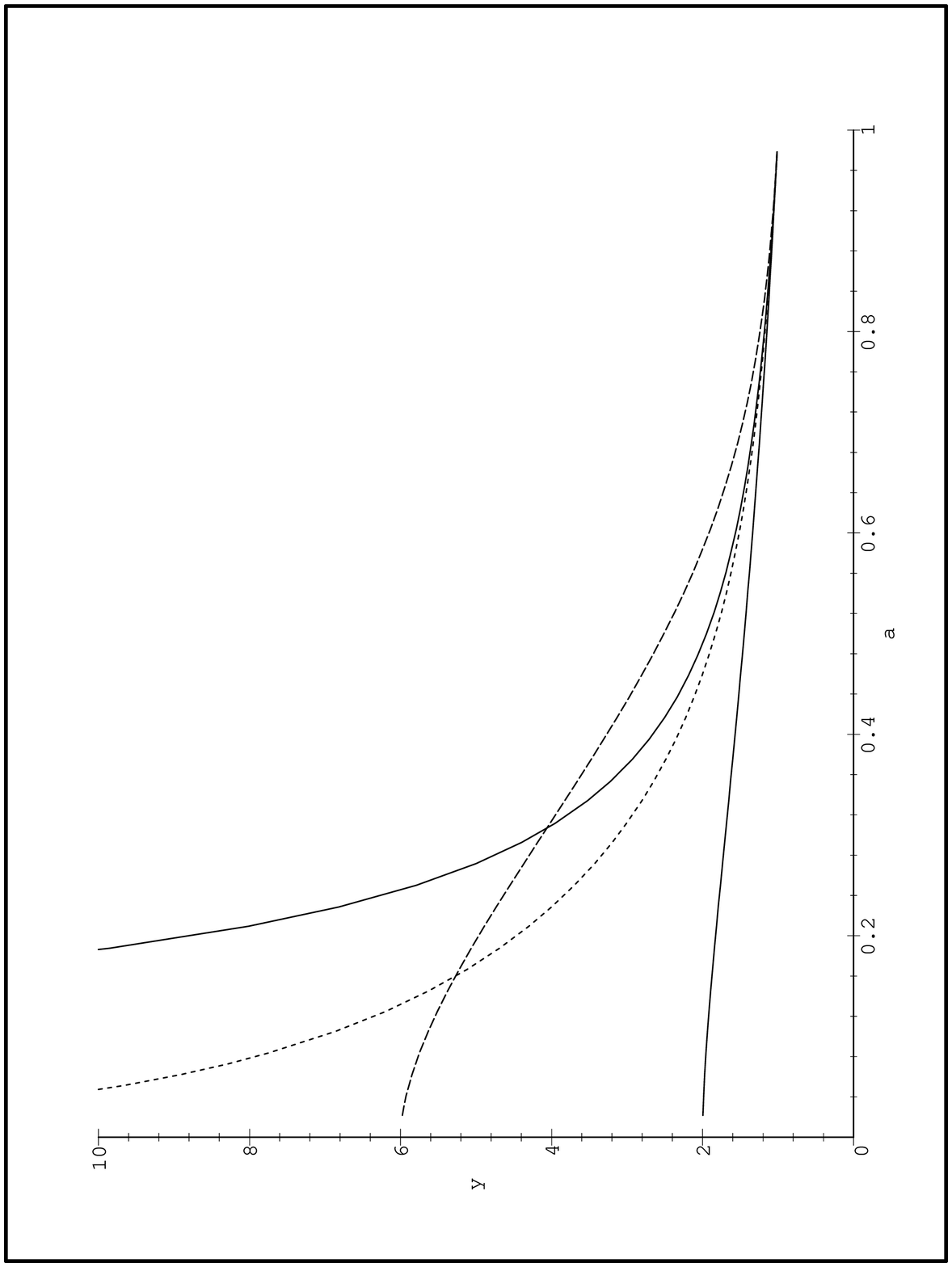,width=10cm}}
\end{rotate}
\]

\vspace{-1cm}

Figure  1 ($0 \le a \le 1$) 
\end{center}

\[
\begin{rotate}[r]{\epsfig{figure=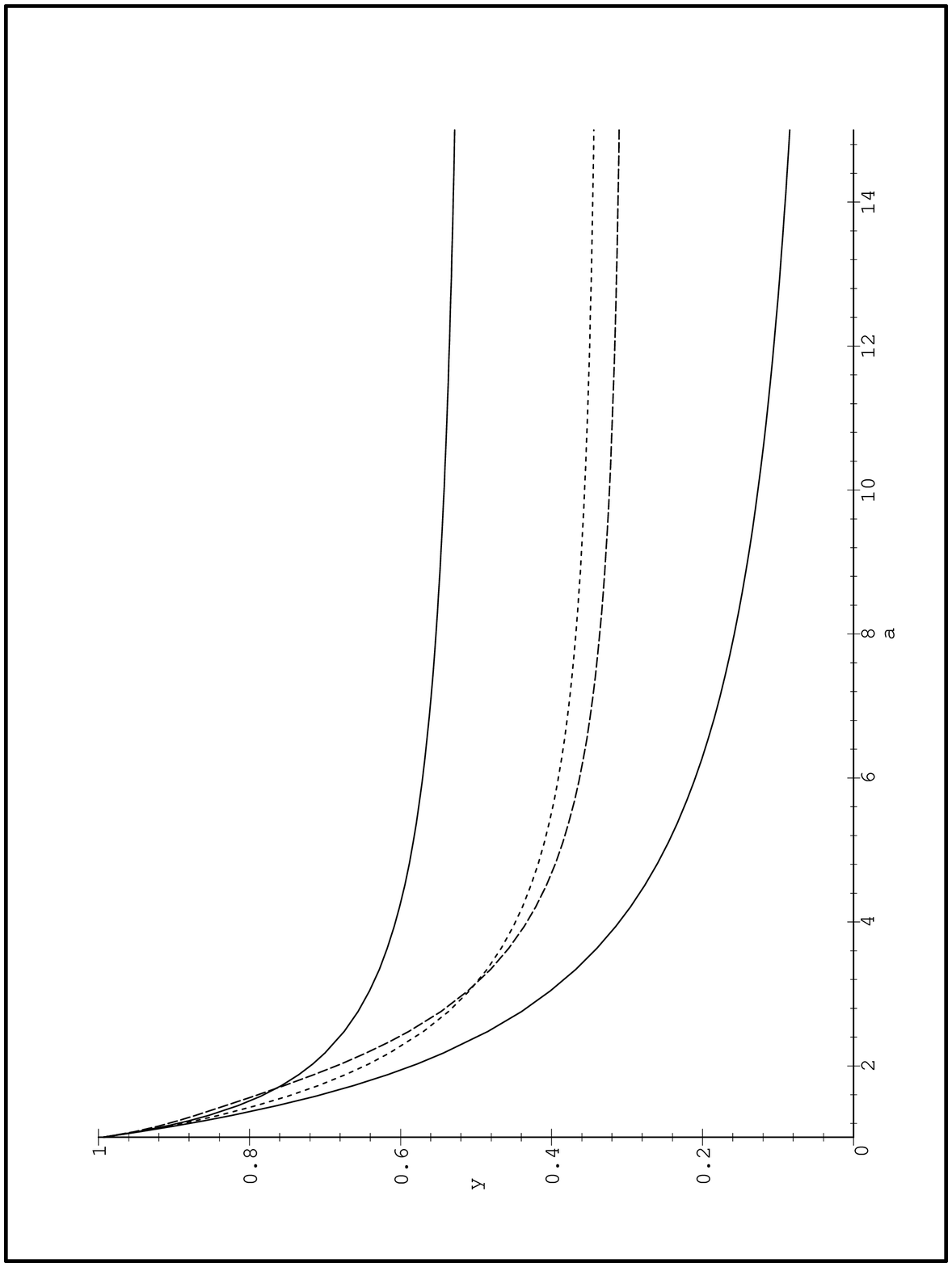,width=10cm}}
\end{rotate}\]

\vspace{-1.5cm}

\begin{center} Figure 2 ($1 \le a \le \infty$) \end{center}

\subsection{Evaluation of the intrinsic upper bound}

We now apply Theorem 2 to the ellipsoid $E(a)$. We can find a uniformization map $\Phi : S^2 \to E(a)$ of the form $ \Phi (x, \varphi)= (w (x), \varphi)$. By formula 1.) for $ds_a^2$ we obtain

\[ \Phi^* (ds^2_a)= \frac{\Delta_a (w(x))}{1-w^2(x)} [w' (x)]^2 dx^2 + (1-w^2 (x)) d \varphi^2 \]

and the condition

\[ \Phi^* (ds^2_a)= h_a^4 (x) \frac{4}{(1+x^2)^2}  \{ dx^2 + x^2d \varphi^2 \} \]

implies the differential equation\\

\mbox{} \hfill $\displaystyle{\frac{\Delta_a^{1/2} (w(x))}{1-w^2} w' = - \frac{1}{x}}$ \hfill $(*)$\\

as well as the boundary conditions $w(0) =1 $ and $ w(\infty)= - 1$. The function $h_a^4 (x)$ is then given by 

\[ h_a^4 (x) = (1- w_a^2 (x)) \frac{(1+x^2)^2}{4x^2}  , \]

where $w_a (x)$ is the unique solution of the differential equation $(*)$ depending on the parameter $a$. We calculate the gradient of $h_a (x)$ with respect to the standard metric

\[ g_o = \frac{4}{(1+x^2)^2} \{ dx^2 + x^2 d \varphi^2 \}  \]

of the sphere $S^2$ and finally obtain

\[ I_1 {(a)} := \int\limits_{S^2} \frac{|\mbox{grad} (h_a)|^2}{h_a^2} dS^2 = \frac{\pi}{2} \int\limits^{\infty}_0 \frac{1}{x} \left( \frac{w_a (x)}{\Delta_a^{1/2} (w_a(x))} + \frac{x^2 -1}{x^2 +1} \right)^2 dx . \]

Theorem 2 then provides the inequality

\[ \lambda_1^2 (a) \le \frac{4 \pi}{\mbox{vol} (E(a))} + \frac{I_1 (a)}{\mbox{vol} (E(a))} . \]

The solution of the differential equation $(*)$ has the symmetry $w_a (x) = - w_a \left( \frac{1}{x} \right)$. Indeed, suppose that  $w_a (x)$ is a solution and consider $w^* (x) = - w_a \left( \frac{1}{x} \right)$. Then $w^*$ solves again the differential equation $(*)$ and $w^* (0) = - w_a (\infty)=1$, \linebreak $w^* (\infty) = - w_a (0) = -1$. This implies that,  for any parameter $0 < a < \infty$,  the solution $w_a (x)$ of the equation $(*)$ vanishes at $x=1$. Consequently $w_a (x)$ is a decreasing function and  we have

\[ \left\{ 
\begin{array}{lll} w_a (x) \ge 0 & \mbox{for} \quad 0 \le x \le 1, & 0 < a < \infty , \\
\mbox{}\\
w_a (x) \le 0 & \mbox{for} \quad 1 \le x \le \infty,  & 0 < a < \infty . \end{array} \right. \]

In particular, $I_1 (a)$ may be reduced to an integral over the interval $[0,1]$:  

\[ I_1 (a) = \pi \int\limits^1_0 \frac{1}{x} \left( \frac{w_a (x)}{\Delta_a^{1/2} (w_a (x))} + \frac{x^2-1}{x^2+1} \right)^2 dx . \]

We study again the limits for $a \to 0, \infty$. First we consider the case that $a \le 1$. Then,  for all points $0 \le x \le 1$,  we have 

\[ \Delta_a^{1/2} (w_a (x)) = \sqrt{(1-a^2) w_a^2 (x) + a^2} \ge \sqrt{1-a^2} \, \, w_a (x) \]

and, consequently, 

\[ \frac{\Delta_a^{1/2} (w_a (x))}{1-w_a^2 (x)} \, \, w'_a (x) \le \sqrt{1-a^2} \, \, \frac{w_a (x) w'_a (x)}{1-w_a^2}  . \]

We integrate this inequality on the interval $[y,1]$. Using the fact that $w_a (1)=0$,  we obtain the estimate

\[ w_a^2 (y) \le 1 - y^{\frac{2}{\sqrt{1-a^2}}} \quad , \quad 0 \le y \le 1 . \]

On the other hand, we have $\Delta_a^{1/2} (w_a (x)) \le 1$. This inequality implies

\[ \frac{w'_a (x)}{1-w^2_a (x)} \le \frac{\Delta_a^{1/2} (w_a (x)) w'_a (x)}{1-w^2_a (x)} = - \frac{1}{x} \]

and, finally, 

\[ w_a (y) \ge \frac{1-y^2}{1+y^2} \quad , \quad 0 \le y \le 1 . \]

Altogether, for any $x \in [0,1]$,  we obtain the inequalities

\[ \frac{1-x^2}{1+x^2} \le \lim\limits_{\overline{a \to 0}}  \, \, w_a (x) \le \overline{\lim\limits_{a \to 0}} \, \, w_a (x) \le 1 -x^2 . \]

Now we apply the following observation: Let $w_a$ be a sequence of numbers such that

\begin{itemize}
\item[a.)] $ 0 < w_a < 1$;
\item[b.)] $ \lim\limits_{\overline{a \to 0}} \, w_a >0$ . 
\end{itemize}

Then the sequence ${w_a}/{\Delta_a^{1/2}}$ with $\Delta_a =( 1-a^2) w_a^2 + a^2$ converges to 1, i.e., 

\[ \lim\limits_{a \to 0} \frac{w_a}{\Delta_a^{1/2}} =1 . \]

In our situation we can conclude that

\[ \lim\limits_{a \to 0} \frac{w_a (x)}{\Delta_a^{1/2} (w_a (x))} = 1  \]

and finally we are able to calculate the limit:

\begin{eqnarray*}
 \lim\limits_{a \to 0} I_1 (a) &=& \lim\limits_{a \to 0} \pi \int\limits^1_0 \frac{1}{x} \left( \frac{w_a (x)}{\Delta^{1/2}_a (w_a (x))} + \frac{x^2-1}{x^2+1} \right)^2 dx = \\
&=&  \pi \int\limits^1_0 \frac{1}{x} \left( 1+ \frac{x^2-1}{x^2+1} \right)^2 dx = 4 \pi \int\limits^1_0 \frac{x^3}{(1+x^2)^2} dx . 
\end{eqnarray*}

Using $\lim\limits_{a \to 0} \mbox{vol} \, (E(a))= 2 \pi$ we obtain

\[ \overline{\lim\limits_{a \to 0}} \lambda_1^2 (a) \le 2+2 \int\limits^1_0 \frac{x^3}{(1+x^2)^2} dx  = \frac{3}{2} + \ln 2 \approx 2,2 . \]

In a similar way we handle the case that $a \ge 1$. The inequalities $(0 \le x \le 1)$

\[ 1 \le \Delta_a (w_a (x)) \le a^2 \]

allow us to prove the estimate

\[ \frac{1-x^{2/a}}{1+x^{2/a}} \le w_a (x) \le \frac{1-x^2}{1+x^2}  , \]

which is valid for all $0 \le x \le 1$ and $a \ge 1$. However,  the function ${w}/{\Delta_a^{1/2} (w)}$ is a monotone decreasing function for $w >0$. Consequently, we have

\[ \D \frac{\frac{1-x^{2/a}}{1+x^{2/a}}}{\Delta_a^{1/2} \left( \frac{1-x^{2/a}}{1+x^{2/a}} \right)} \le \frac{w_a (x)}{\Delta_a^{1/2} (w_a (x))} \le 1  \]

and from this inequality we can deduce 

\[ \left( 1 - \frac{w_a(x)}{\Delta_a^{1/2} (w_a (x))} \right)^2 \le \frac{4a^2 x^{2/a}}{(1-a^2)(1-x^{2/a})^2 + a^2 (1+x^{2/a})^2} . \]

We split the integral $I_1 (a)$ into three parts:

\begin{eqnarray*}
I_1 (a) &=& \pi \int\limits_0^1 \left( \left( \frac{w_a (x)}{\Delta_a^{1/2} (w_a (x))} - 1 \right) + \frac{2x^2}{1+x^2} \right)^2 dx =\\
&=& 4 \pi \int\limits^1_0 \frac{x^3}{(1+x^2)^2}dx + 4 \pi \int\limits^1_0 \left( \frac{w_a (x)}{\Delta_a^{1/2} (w_a (x))} - 1 \right) x dx \\
&& + \pi \int\limits^1_0 \frac{1}{x} \left( \frac{w_a (x)}{\Delta_a^{1/2} (w_a (x))} - 1 \right)^2 dx . 
\end{eqnarray*}

We estimate the last term using the inequality for $\left( 1 - \frac{w_a (x)}{\Delta_a^{1/2} (w_a (x))} \right)^2$ and obtain as the value of the integral

\[ I_1 (a) \le 4 \pi \int\limits^1_0 \frac{x^3}{(1+x^2)^2} dx + 4 \pi \int\limits^1_0 \left( \frac{w_a (x)}{\Delta^{1/2}_a (w_a (x))}  -1 \right) xdx +  \]
\[  +  \frac{\pi}{2} a \cdot \frac{a}{\sqrt{a^2 - 1}} \ln \left( \frac{8a^3 - 6a^2 + 2a \sqrt{a^2 - 1}}{8a^3 - 6a^2 - 2a \sqrt{a^2 - 1}} \right) . \]

The volume $\mbox{vol} \, (E(a))$ of the ellipsoid behaves like $\pi^2a$, i.e., 

\[ \lim\limits_{a \to \infty} \frac{\mbox{vol} \, (E(a))}{a} = \pi^2 . \]

Therefore,  we can control the asymptotic behaviour  of $\lambda_1^2 (a)$ for $a \to \infty$:

\[ \lambda_1^2 (a)  \le  \frac{4 \pi}{a} \frac{a}{\mbox{vol} \, (E(a))} + \frac{I_1 (a)}{a} \frac{a}{\mbox{vol} \, (E(a))} \stackrel{\le}{\sim}  \frac{4}{\pi  a} + \frac{4}{\pi  a} \int\limits^1_0 \frac{x^3}{(1+x^2)^2} dx = \frac{1}{\pi a} [ 2 \ln (2) + 3 ] . 
\]

In particular,  we have shown

\[ \lim\limits_{a \to \infty} \lambda_1^2 (a) = 0 . \]

\section{The first eigenvalue of the Dirac operator on the tube around a circle}

We consider a circle in a plane with curvature $\kappa$ and length $L = 2 \pi / \kappa$. Let $r$ be a fixed radius and denote by $M^2 (r)$ its tube in ${\Bbb R}^3$ of radius $r$, $r \kappa <1$. The induced metric on the surface $M^2 (r)$ is given by the formula

\[ g = (1-r \kappa \cos \varphi)^2 ds^2 + r^2 d \varphi^2  , \]

where we use the length parameter $0 \le s \le L$ for the circle and $0 \le \varphi \le 2 \pi$ parametrizes the angle of the tube. First of all we calculate a uniformization

\[  \Phi : [0,L] \times [0,A] \to [0,L] \times [0, 2 \pi] \]

of this metric on $T^2$. Suppose $\Phi$ is given by the condition $\Phi (s, \psi) =(s, \varphi (\psi))$. Then the equation $\Phi^* (g)= h^4 (ds^2 + d \psi^2)$ yields the differential equation 

\[ \frac{\varphi' (\psi)}{1-r \kappa \cos (\varphi (\psi))} = \frac{1}{r} \]

and the function $h = h(s, \psi)$ is given by

\[ h^2 = r \varphi' (\psi) = 1- r \kappa \cos (\varphi (\psi)) . \]

Using the integral $(a<1)$

\[ \int \frac{dx}{1-a \cos (x)} = \frac{2}{\sqrt{1-a^2}} \, \, \mbox{arc} \,  \mbox{tg} \,  \left( \frac{ (1+a) \mbox{tg} \,  \left( \frac{x}{2} \right)}{\sqrt{1-a^2}} \right) \]

we obtain the solution $\varphi (\psi)$

\[ \mbox{tg} \, \left( \frac{\varphi (\psi)}{2} \right) = \sqrt{ \frac{1-r \kappa}{1+ r \kappa}} \, \, \mbox{tg} \, \left( \frac{1}{2r} \sqrt{1-r^2 \kappa^2} \, \psi \right) . \]

Since $\varphi (\psi)$ maps the interval $[0,A]$ bijectively onto $[0,2 \pi]$,  we conclude \linebreak $A= 2\pi r / \sqrt{1-r^2 \kappa^2}$. Moreover, the function $h^2$ is determined by

\begin{eqnarray*}
h^2 &=& 1-r \kappa \cos (\varphi (\psi)) = 1 - r \kappa \, \, \frac{1- \mbox{tg}^2 \, \left(  \frac{\varphi (\psi)}{2} \right)}{1+ \mbox{tg}^2 \, \left(  \frac{\varphi (\psi)}{2} \right)} =\\
\mbox{}\\
&=& (1-r^2 \kappa^2) \, \, \frac{1 + \mbox{tg}^2 \, \left(\frac{1}{2r} \sqrt{1-r^2 \kappa^2} \, \psi \right)}{1+r \kappa + (1-r \kappa) \mbox{tg}^2 \, \left(\frac{1}{2r} \sqrt{1-r^2 \kappa^2}\,  \psi \right)} . 
\end{eqnarray*}

Hence, we obtain a uniformization of the metric of the tube $M^2 (r)$ parametrized on $[0,L] \times \left[0, \frac{2 \pi r}{\sqrt{1-r^2 \kappa^2}} \right]$. The basis of the lattice is

\[ v_1 = \left( L , 0  \right) \quad , \quad v_2 = \left(  0 ,   \frac{2 \pi r}{\sqrt{1-r^2 \kappa^2}}  \right) , \]

and thus the dual lattice has the basis

\[  v_1^* = \left(  \frac{1}{L}  , 0  \right) \quad , \quad v_2^* = \left(  0 ,   \D \frac{\sqrt{1-r^2 \kappa^2}}{2 \pi r}  \right) . \]

By Theorem 4 we obtain the estimate

\[ \lambda_1^2 (\kappa, r) \le \frac{1}{4} \left( \kappa^2 \varepsilon_1 + \frac{1-r^2 \kappa^2}{r^2} \varepsilon_2 \right) \, \, \frac{\D \int_{T^2 \downsp} \frac{1}{h^2} dT^2}{\D \int_{T^2}^{\upsp} h^2 dT^2} \]

for the first eigenvalue of the Dirac operator on the tube $M^2(r)$ with respect to the spin structure  $(\varepsilon_1, \varepsilon_2)$. We compute these two integrals:

\[ \int_{T^2} h^2 dT^2 = \int\limits^L_0 \int\limits^A_0 h^2 (s, \psi) ds d \psi =  L \cdot \int\limits^A_0 r \varphi' (\psi) d \psi = L r \int\limits^{2 \pi}_0 d \varphi = 2 \pi r L , 
\]

and

\[
\int_{T^2} \frac{1}{h^2} d T^2 = \frac{L}{r} \int\limits^A_0 \frac{1}{\varphi' (\psi)} d \psi  = \frac{L}{r} \int\limits^A_0 \frac{r^2}{(1-r \kappa \cos (\varphi (\psi)))^2} \varphi' (\psi) d \psi =  Lr \int\limits^{2 \pi}_0 \frac{d \varphi}{(1-r \kappa \cos (\varphi))^2} \quad . \]

Consequently, this ratio is equal to

\[ \frac{\D \int_{T^2 \downsp} \frac{1}{h^2} \, \, dT^2}{\D \int_{T^2}^{\upsp} h^2 dT^2} = \frac{1}{2 \pi} \int\limits^{2 \pi}_0 \frac{d \varphi}{(1- r \kappa \cos (\varphi))^2} = \frac{1}{(1 - r^2 \kappa^2)^{3/2}} \, \, , \]

i.e., 

\[ \lambda_1^2 (\kappa, r) \le \frac{1}{4} \left( \kappa^2 \varepsilon_1 + \frac{1-r^2 \kappa^2}{r^2} \varepsilon_2  \right) \frac{1}{(1-r^2 \kappa^2)^{3/2} } \, \, . \]

\newcommand{\vol}{\mbox{vol} \, }
\newcommand{\grad}{\mbox{grad} \, }

The volume $\vol (\kappa, r)$ of the tube equals

\[ \vol (\kappa, r) = 4 \pi^2 \frac{r}{\kappa} \]

and we obtain the inequality

\[ \lambda_1^2 (\kappa, r) \vol (\kappa,r) \le \pi^2 \left( r \kappa \varepsilon_1 + \frac{1-r^2 \kappa^2}{r \kappa} \varepsilon_2 \right) \frac{1}{(1-r^2 \kappa^2)^{3/2}} \, \, . \quad \quad (*) \]

Now we apply the inequality

\[ \lambda_1^2 (\kappa,r) \vol (\kappa,r) \le \lambda_1^2 (g_o) \vol (T^2, g_o) + \int_{T^2} \frac{|\grad (h)|^2}{h^2} dT^2 \]

to our situation. Since $h^2 = r \varphi' (\psi) = 1- r \kappa \cos (\varphi (\psi))$,  we can calculate the gradient of $h$:

\[ \frac{|\grad (h)|^2}{h^2} = \frac{r \kappa^2}{4} \, \, \frac{\sin^2 (\varphi (\psi)) \varphi' (\psi)}{1-r \kappa \cos (\varphi (\psi))} \]

and,  therefore,  we obtain

\[ \int_{T^2} \frac{|\grad (h)|^2}{h^2} = \frac{\pi}{2} r \kappa \int\limits^{2 \pi}_0 \frac{\sin^2 (\varphi)}{1 - r \kappa \cos (\varphi)} d \varphi = \frac{\pi^2}{r \kappa} \left( 1 - \sqrt{1-r^2 \kappa^2} \right) . \]

Then we have proved the estimate

\[ \lambda_1^2 (\kappa, r) \vol (\kappa,r) \le \pi^2 \left( r \kappa \varepsilon_1 + \frac{1-r^2 \kappa^2}{r \kappa} \varepsilon_2 \right) \cdot \frac{1}{\sqrt{1-r^2 \kappa}} + \frac{\pi^2}{r \kappa} (1 - \sqrt{1-r^2 \kappa^2} ) \, \, . \quad (**) \]

We discuss the inequalities $(*)$ and $(**)$ for the three non-trivial spin structures on the tube. For all cases, we provide a picture in which the long dashed line represents the estimate$(*)$, and the short dashed line the estimate $(**)$. The $x$-axis uses the variable $a=r \kappa$, the $y$-axis is to be understood in multiples of $\pi^2$. For comparison matters only, we have also drawn the line for  constant value 2.\\

{\bf Case 1:} \quad $\varepsilon_1 =1 , \, \, \varepsilon_2 =0$. In this case we obtain

\[ \begin{array}{ccc}  \lambda_1^2 (\kappa,r ) \vol (\kappa,r ) \le \pi^2 r \kappa \D \frac{1}{(1-r^2 \kappa)^{3/2}} & &  (*) \\
\mbox{}\\
 \lambda_1^2 (\kappa, r) \vol (\kappa, r) \le \D \frac{\pi^2 r \kappa}{\sqrt{1-r^2 \kappa^2}} + \frac{\pi^2}{r \kappa} \left(1- \sqrt{1-r^2 \kappa^2}\right) & &  (**)  \end{array} \]

In particular,  we conclude

\[ \lim\limits_{r \to 0} \lambda_1^2 (\kappa,r) \vol (\kappa,r) = \lim\limits_{\kappa \to 0} \lambda_1^2 (\kappa,r) \vol (\kappa,r)=0 . \]

\[
\begin{rotate}[r]{\epsfig{figure=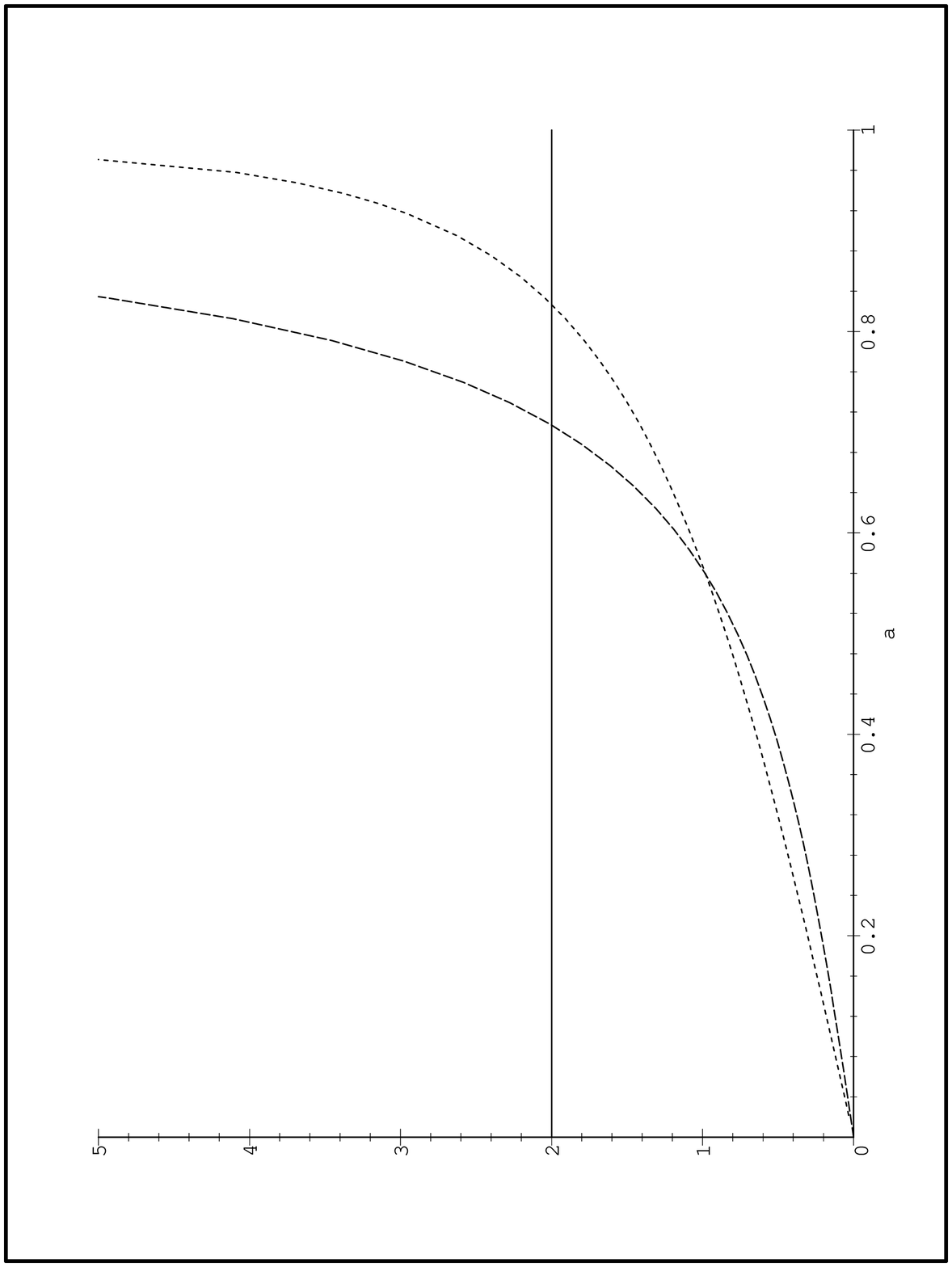,width=10cm}}
\end{rotate} \]

\vspace{-1.5cm}

\begin{center} Figure 3 ($\varepsilon_1 =1, \varepsilon_2 =0$)  \end{center}

\newpage

{\bf Case 2:} \quad $\varepsilon_1 =0 , \, \, \varepsilon_2 =1$. In this case the inequalities are

\[ \begin{array}{ccc} \lambda_1^2 (\kappa,r) \vol (\kappa,r) \le \D \frac{\pi^2}{r \kappa} \frac{1}{\sqrt{1-r^2 \kappa^2}} & & (*) \\
\mbox{}\\
\lambda_1^2 (\kappa, r) \vol (\kappa,r) \le \D \frac{\pi^2}{r \kappa} & & (**)  \end{array} \]

and, in particular,  we conclude 

\[ \overline{\lim\limits_{r \kappa \to 1}} \, \, \lambda_1^2 (\kappa , r) \vol (\kappa, r) \le \pi^2 . \]

\[
\begin{rotate}[r]{\epsfig{figure=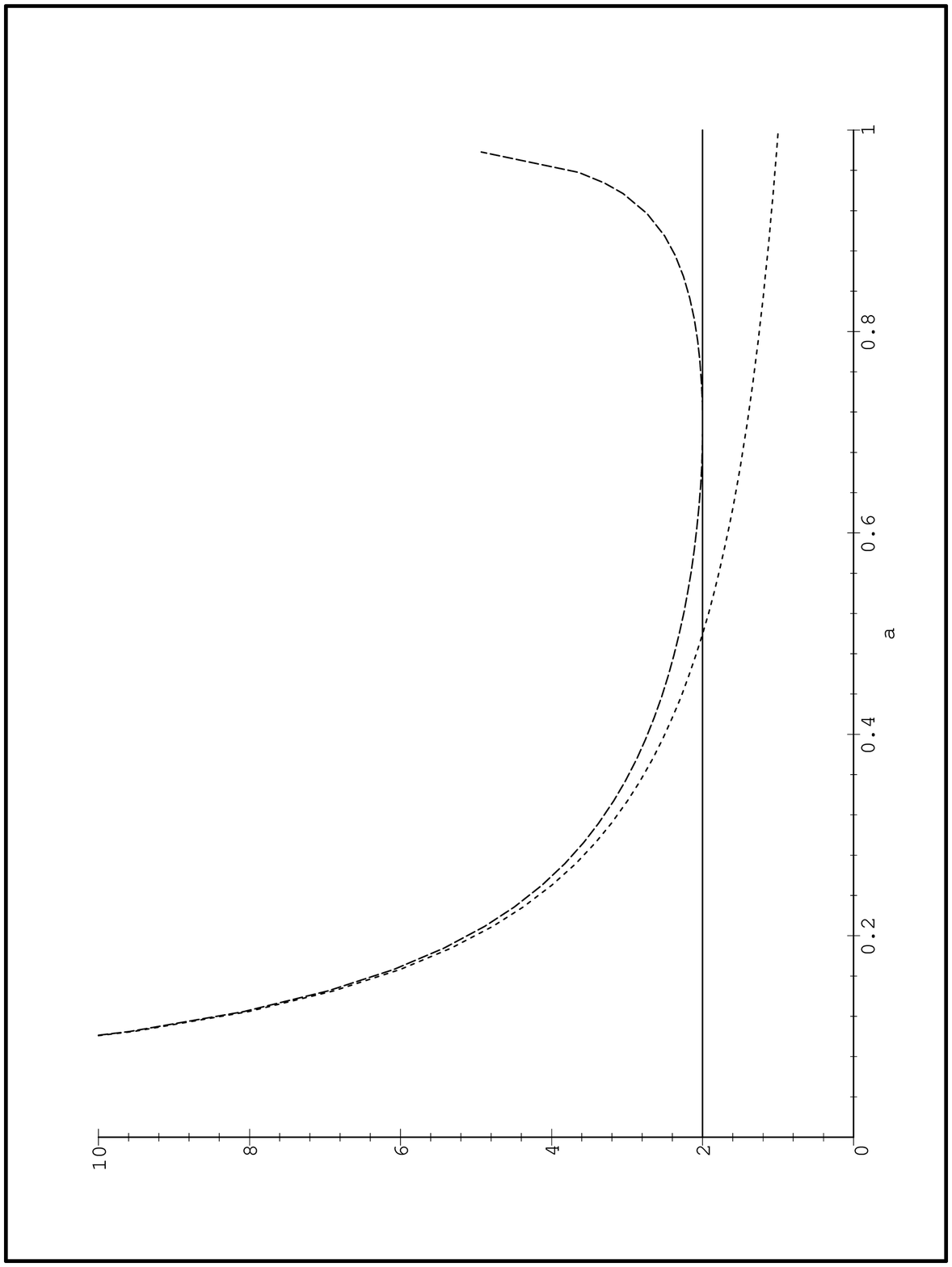,width=10cm}}
\end{rotate}\]

\vspace{-1.5cm}

\begin{center} Figure 4 ($\varepsilon_1=0 , \varepsilon_2 =1$) \end{center}

\vspace{2cm}

{\bf Case 3:} \quad $\varepsilon_1 = 1 = \varepsilon_2$. \, \, In this case we obtain the estimates

\[ \mbox{} \hspace{1.7cm} \lambda_1^2 (\kappa,r) \vol (\kappa,r) \le \frac{\pi^2}{r \kappa} \frac{1}{(1-r^2 \kappa^2)^{3/2}} \hspace{2.6cm} (*) \]

\[ \lambda_1^2 (\kappa, r) \vol (\kappa, r) \le  \frac{\pi^2}{r \kappa} \frac{1}{\sqrt{1-r^2 \kappa^2}} + \frac{\pi^2}{r \kappa} \left(1 - \sqrt{1-r^2 \kappa^2} \right) . \quad \quad (**) \]

Let us compare these estimates obtained via the uniformization of the tube with the estimate using the embedding $M^2 (r) \subset {\Bbb R}^3$. Notice that the embedding induces the spin structure $\varepsilon_1 = 1 = \varepsilon_2$ on the tube. Then we obtain

\[ \lambda_1^2 (\kappa,r) \vol (\kappa,r) \le \int_{M^2(r)} H^2 d M^2 (r) \le \frac{\pi^2}{r \kappa} \frac{1}{\sqrt{1-r^2 \kappa^2}} \, , \quad \quad (***) \]

i.e.,  the extrinsic bound (drawn as a solid line in figure 5) for $\lambda_1^2$ \, is better than the intrinsic estimates.\\

\[
\begin{rotate}[r]{\epsfig{figure=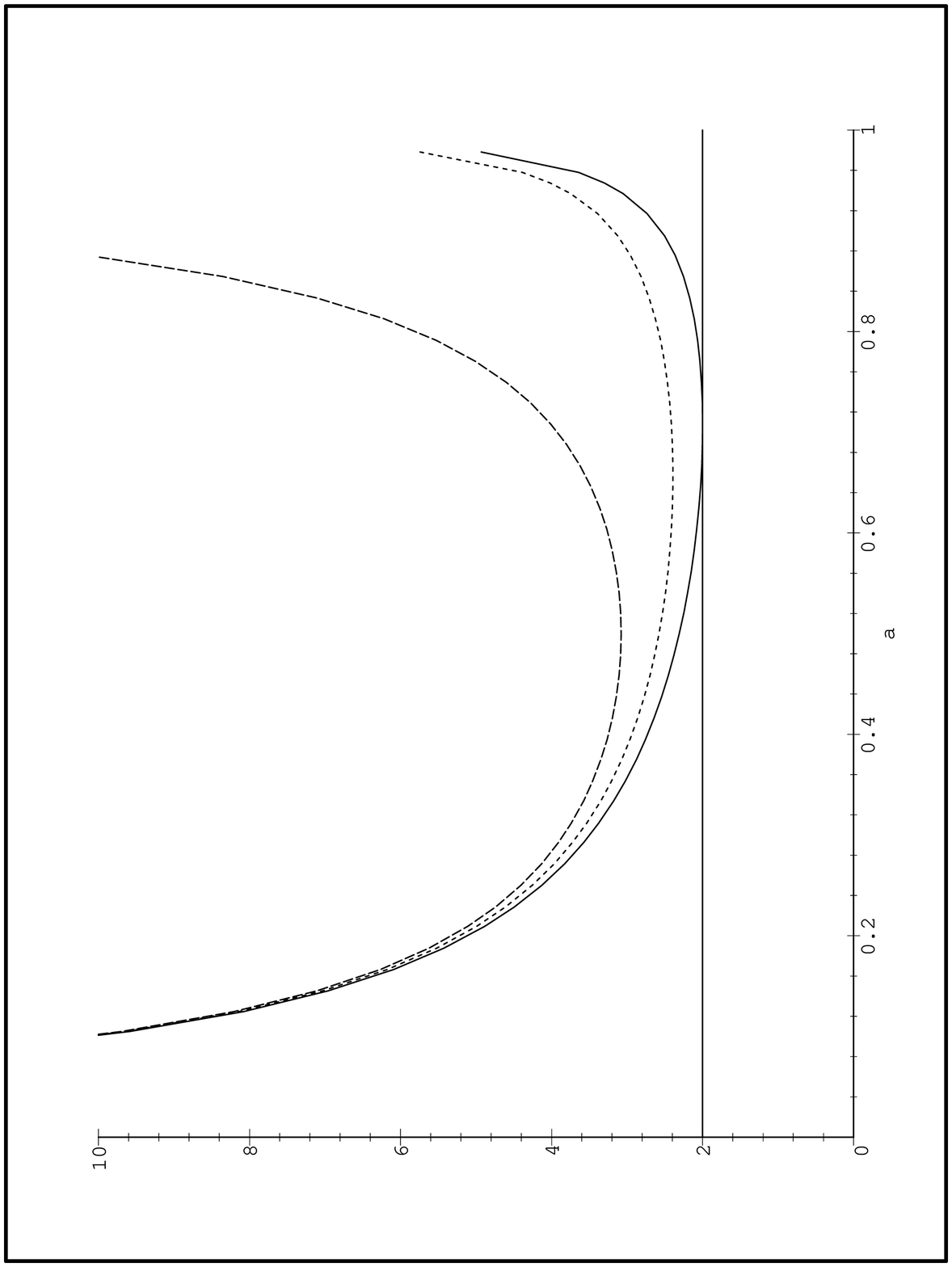,width=10cm}}
\end{rotate}\]

\vspace{-1.5cm}

\begin{center} Picture 5 ($\varepsilon_1 = 1 = \varepsilon_2$) \end{center}

\vspace{1cm}

{\bf Case 4:} \quad $\varepsilon_1 = 0 = \varepsilon_2$ \, \, \, In this case $\lambda_0 (D)=0$ is an eigenvalue of the Dirac operator and Theorem 3 yields the following estimate for the first non-trivial eigenvalue $\lambda_1^2 (\kappa,r)$:

\[ \lambda_1^2 (\kappa,r) \le \min \left\{ 2 \kappa^2 , \frac{1}{r^2} \right\} \, \,  \frac{\D \int\limits^{2 \pi}_{0 \downsp} \frac{d \varphi}{(1-r \kappa \cos \varphi)^4}}{\D \int\limits^{2 \pi^{\upsp} }_0 \frac{d \varphi}{(1-r \kappa \cos \varphi)^2}} . \]

In particular,  we obtain

\[ \lim\limits_{\kappa \to 0} \lambda_1^2 (\kappa,r)=0 \quad , \quad \overline{\lim\limits_{r \to 0}} \, \, \lambda_1^2 (\kappa,r) \le 2 \kappa^2 \]

and

\[ \lim\limits_{r \cdot \kappa \to 0} \lambda_1^2 (\kappa,r) \vol (\kappa,r)=0 . \]


\end{document}